\RequirePackage{silence}
\documentclass[10pt, english]{article}

\usepackage[margin= 2 cm]{geometry}

\usepackage{amsthm}
\usepackage{amsmath}
\usepackage{amssymb}
\usepackage{setspace}
\usepackage{mathtools}
\usepackage{verbatim}
\usepackage{csquotes}
\usepackage{graphicx}
\usepackage[hidelinks]{hyperref}
\usepackage{thm-restate}
\usepackage{cleveref}
\usepackage{graphicx}
\usepackage{enumitem}
\usepackage{framed}
\usepackage{subcaption}

\usepackage{floatrow}
\usepackage[T1]{fontenc}
\usepackage{tikz}
\usepackage{mathdots}
\usepackage{xcolor}
\usepackage{diagbox}
\usepackage{colortbl}
\usepackage[absolute,overlay]{textpos}

\graphicspath{ {./images/} }

\usetikzlibrary{calc}
\usetikzlibrary{decorations.pathreplacing}
\usetikzlibrary{positioning,patterns}
\usetikzlibrary{arrows,shapes,positioning}
\usetikzlibrary{decorations.markings}

\tikzstyle{edge}=[very thick]
\definecolor{bostonuniversityred}{rgb}{0.8, 0.0, 0.0}
\definecolor{arsenic}{rgb}{0.23, 0.27, 0.29}
\tikzstyle{diredge}=[postaction={decorate,decoration={markings,
		mark=at position .95 with {\arrow[scale = 1]{stealth};}}}]

\newcommand{\fitellipsis}[2] 
{\draw [fill=green]let \p1=(#1), \p2=(#2), \n1={atan2(\y2-\y1,\x2-\x1)}, \n2={veclen(\y2-\y1,\x2-\x1)}
    in ($ (\p1)!0.5!(\p2) $) ellipse [ x radius=\n2/2+0cm, y radius=0.1cm, rotate=\n1];
}
\newcommand{\cali}[1]{\mathcal{#1}}

\theoremstyle{plain}

\newtheorem*{thm*}{Theorem}
\newtheorem{thm}{Theorem}
\Crefname{thm}{Theorem}{Theorems}
\numberwithin{thm}{subsection}

\newtheorem*{lem*}{Lemma}

\Crefname{lem}{Lemma}{Lemmas}

\newtheorem*{claim*}{Claim}

\Crefname{claim}{Claim}{Claims}
\Crefname{claim}{Claim}{Claims}

\newtheorem{prop}[thm]{Proposition}
\Crefname{prop}{Proposition}{Propositions}

\newtheorem{cor}[thm]{Corollary}
\Crefname{cor}{Corollary}{Corollaries}

\newtheorem{conj}[thm]{Conjecture}
\Crefname{conj}{Conjecture}{Conjectures}

\Crefname{qn}{Question}{Questions}

\newtheorem{obs}[thm]{Observation}
\Crefname{obs}{Observation}{Observations}

\newtheorem{ex}[thm]{Example}
\Crefname{ex}{Example}{Examples}

\theoremstyle{definition}

\Crefname{prob}{Problem}{Problems}

\Crefname{defn}{Definition}{Definitions}

\newtheorem*{defn*}{Definition}

\theoremstyle{remark}

\renewenvironment{proof}[1][]{\begin{trivlist}
\item[\hspace{\labelsep}{\bf\noindent Proof#1.\/}] }{\qed\end{trivlist}}

\expandafter\def\expandafter\normalsize\expandafter{%
    \normalsize
    \setlength\abovedisplayskip{8pt}
    \setlength\belowdisplayskip{8pt}
    \setlength\abovedisplayshortskip{4pt}
    \setlength\belowdisplayshortskip{4pt}
}

\usepackage[square,sort,comma,numbers]{natbib}
 \setlist[itemize]{leftmargin=*}

\newcommand{\T}[1]{\texttt{#1}}

\newcommand{\fitellipsiss}[3] 
{\draw [fill=white]let \p1=(#1), \p2=(#2), \n1={atan2(\y2-\y1,\x2-\x1)}, \n2={veclen(\y2-\y1,\x2-\x1)}
    in ($ (\p1)!0.5!(\p2) $) ellipse [ x radius=\n2/2+#3cm, y radius=#3cm, rotate=\n1];
}

\DeclareFontFamily{OT1}{pzc}{}
\DeclareFontShape{OT1}{pzc}{m}{it}{<-> s * [1.10] pzcmi7t}{}

\title{\vspace{-0.8cm} Automated  Counting of Spanning Trees for  Several Infinite families of Graphs}
\author{Pablo Blanco \and Doron Zeilberger}
 \date{}

\begin{document}

\usetikzlibrary{decorations.pathreplacing}
\usetikzlibrary{positioning,patterns}
\usetikzlibrary{arrows,shapes,positioning}
\usetikzlibrary{decorations.markings}

\maketitle

\vspace{-0.5cm}
\begin{abstract}
    Using the theoretical basis developed by Yao and Zeilberger, we consider certain graph families whose structure results in a rational generating function for sequences related to spanning tree enumeration. 
Said families are powers of cycles and powers of paths; later, we briefly discuss toru graphs and grid graphs.
 In each case we know, a priori, that the 
set of spanning trees of the family of graphs can be described in terms of a finite-state-machine, and hence
there is a finite transfer-matrix that guarantees the generating function is rational. Finding this ``grammar'', and hence the transfer-matrix is very tedious, so a much more efficient
approach is to use experimental mathematics. Since computing numerical determinants is so fast, one can use the matrix tree theorem to generate sufficiently many terms, then fit the data to a rational function. The whole procedure can be done rigorously \textit{a posteriori}.

We also construct generating functions for the quantity ``total number of leaves'' over all spanning trees, and automatically derive the asymptotics for both number of spanning trees and
the sum of the number of leaves, enabling our computer to get the asymptotic for the average number of leaves per vertex in a random spanning tree in each family that always tends
to a certain number, which we compute explicitly. This number, that we christen the B-Z constant for the infinite family is always some (often complicated) algebraic number, in contrast
to the family of complete graphs where this constant is $1/e$ (a transcendental number).
    \end{abstract}

\section{Introduction}
A subgraph $T$ of a graph $G$ such that $T$ is a tree with $V(T)=V(G)$ is called a \textit{spanning tree} of $G$. If $G$ is connected, simple, and not a tree, then it contains several spanning trees. When we count the spanning trees of a graph in this paper, we consider two isomorphic trees with different vertex labels to be different trees; for example, the complete graph on three vertices $K_3$ has three spanning trees, not one. For a connected graph $G$, we denote its number of spanning trees by $\tau(G)$. Cayley's formula counts the number of spanning trees of a complete graph on $n$ vertices to be $n^{n-2} = \tau(K_n)$.

The \textit{generating function} $f(x)$ for a sequence $(a_0,a_1,\dots)$ is the formal power series obtained by setting the sequence terms as its coefficients:
$$f(x):=\sum_{n=0}^\infty a_nx^n.$$
Whenever a sequence has a recurrence of finite length, with constant coefficients, it is called \textit{C-finite}. In other words, $(a_0,a_1,\dots)$ is C-finite (of order $r$) if there is a fixed $r$ and coefficients $c_0,\dots,c_{r-1}$ with $c_0\neq0$ such that
$$a_{n+r}=c_{r-1}a_{n+r-1} + \dots +c_0a_n $$
for all $n\ge 0$.

The following is a useful property of C-finite sequences:
\begin{thm}
    \emph{\cite[Thm 4.3]{concrete}} A sequence $(a_0,a_1,\dots)$ is C-finite (of order $r$) with recurrence
    $$a_{n+r}=c_{r-1}a_{n+r-1} + \dots +c_0a_n $$
    if and only if
    $$\sum a_nx^n = \frac{p(x)}{1-c_{r-1}x+\dots-c_1x^{r-1}-c_0x^r}$$
    for some polynomial $p(x)$ of degree at most $r-1$. 
\end{thm}
 In fact, the polynomial $p(x)$ in the previous theorem is determined by the initial values $(a_0,\dots,a_{r-1})$ of the sequence. 
Later, we will introduce various (ordered) families of graphs. Each graph family will give rise to several numerical sequences (such as number of spanning trees) which will be indexed, in some way, by the number of vertices of the underlying graph in the family. 
To disambiguate the numerator $p(x)$, we will specify the smallest member in each graph family (which corresponds to the first term, $a_0$, in any resulting numerical sequence).

 Due to the recursive nature of the graph families we consider, we presume that there is a finite transfer matrix which describes the corresponding sequence. (In a paper by Yao and Zeilberger \cite{gordianknot}, there is an example of such a construction for the family of $k\times n$ grid graphs, for fixed $k$.) Thanks to this, we can find the sequence's rational generating function by computing a large number of terms in the sequence. The guessing is done in a Maple procedure, described in the same paper \cite{gordianknot}, called \texttt{GuessRec}.

 To generate the terms of whichever sequence we consider, we apply Kirchhoff's matrix-tree theorem, which allows us to compute the number of spanning trees of a(n $n$-vertex) graph by starting with its Laplacian matrix, taking an $n-1$ by $n-1$ (matrix) minor, and computing its determinant. \Cref{sec:CountingLeaves} provides more detail on how to count leaves by applying matrix-tree theorem. Below, we define the Laplacian matrix of a graph for the reader's convenience.

 If $G$ is a (simple) graph with vertices $v_1,\dots,v_n$, its \emph{Laplacian matrix} is a symmetric matrix with its entries defined by
 \begin{equation*}
     a_{i,j}:=\begin{cases}
         -1 & \text{if } v_i\sim v_j \text{ and } i\neq j\\
    \deg(v_i) & \text{if } i=j \\
     \end{cases}.
 \end{equation*}

 If the reader wishes to further review labelled spanning trees, we direct them to classics such as \cite{aldous} or \cite{moon}.

In \Cref{sec:CountingLeaves}, we describe methods to count the leaves (across all spanning trees) for graphs in the same families. For \Cref{sec:bz-constant}, we consider an asymptotic constant relating total number of leaves to total number of spanning trees in any families (for which the constant is well-defined). We state our main experimental results in \Cref{sec:cycles} and \Cref{sec:paths}, then briefly consider an additional pair of graph families in \Cref{sec:gridtorus}. In \Cref{sec:verification}, we describe how to numerically estimate the aforementioned constant with random sampling and an implementation to make such estimates. We developed a Maple package that carries out the computations in this paper. It is described in
\Cref{maple-package} and available at \cite{spt}. Finally, we conclude with some conjectures in \Cref{sec:conj}.

\section{Preliminaries} \label{graph-fams}
For ease of notation, we assume all graphs are simple. We write $[n]:=\{1,\dots,n\}$.

\textbf{Definition.} For a graph $G$, the \textit{distance} $|u,v|$ between two vertices $u,v$ is the length of the shortest path between them. If $u$ and $v$ are in different components of $G$, we say the distance between them is infinite and write $|u,v|=\infty$.

\begin{figure}[H] 
   \centering
    \begin{tikzpicture}

\tikzstyle{knode}=[circle,fill,draw=black,thick,inner sep=2pt] 

\pgfmathsetmacro{\n}{7} 
\pgfmathsetmacro{\R}{1cm} 
\pgfmathsetmacro{\A}{90} 

\node (r) at (0,0) {};
\node (s) at (6,0) {};

\foreach \i in {1,...,\n}
{
\node (\i) at ($(r)+(\i*360/\n + \A:1.5cm)$) [knode] {};
}

\pgfmathtruncatemacro{\m}{\n - 1}
\foreach \j in {1,...,\m}
{
\pgfmathtruncatemacro{\k}{\j + 1}
\draw[line width=2pt] (\j)--(\k);
}

\draw[line width=2pt] (\n)--(1);
\draw (1)--(3);
\draw (2)--(4);
\draw (6)--(4);
\draw (3)--(5);
\draw (2)--(7);
\draw (1)--(6);
\draw (7)--(5);

\foreach \i in {8,...,13}
{
\node (\i) at ($(s)+(\i*360/6 + \A:1.5cm)$) [knode] {};
}

\pgfmathtruncatemacro{\m}{12}
\foreach \j in {8,...,\m}
{
\pgfmathtruncatemacro{\k}{\j + 1}
\draw[line width=2pt] (\j)--(\k);
}

\draw (8)--(10);
\draw (8)--(11);
\draw (9)--(11);
\draw (9)--(12);
\draw (10)--(13);
\draw (10)--(12);
\draw (11)--(13);

\end{tikzpicture}
    \caption{On the left, ${C_7}^{2}$. On the right, ${P_6}^{3}$. The thicker edges represent the edges from the corresponding original graph.}
    \label{fig:fig1}
\end{figure}
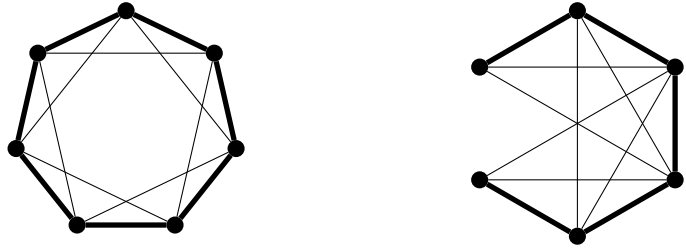

\textbf{Definition.} Let $G$ be a graph. For an integer $k\ge1$, the $k$-th \textit{power} of $G$, denoted $G^k$, is the graph obtained from $G$ such that $V(G^k):=V(G)$ and $E(G^k):= \{uv: u,v\in V(G),\ 1\leq |u,v|\leq k\}$. See \Cref{fig:fig1}.

From the definition, we see that if $d$ is the \emph{diameter} of a connected graph $G$ (defined as $d:=\max_{u,v\in V(G)} |u,v|$ and denoted by $\text{diam }G$), then $G^d$ is the complete graph on $|V(G)|$ vertices. In some cases, it's natural to exclude complete graphs and hence consider only $G^k$ when $k<d$. Additionally, one may note that if there are two shortest paths (in $G$) between two vertices $u,v$ of length $m> 1$, then $G^m$ may be better viewed as a multigraph where there are two edges between $u,v$. 

\textbf{Notation.} We denote the path graph on $n$ vertices by $P_n$ and the cycle graph on $n$ vertices by $C_n$.

In subsequent sections, we state the generating functions for the number of spanning trees in the graph families 
\begin{align*}
    \cali{C}_k&:=\{{C_n}^k\ :\ k < \text{diam } C_n \}\\  \cali{P}_k&:=\{{P_n}^k\ :\ k < \text{diam } P_n \}.
\end{align*}
The condition in each construction ensures that our graph families do not contain any ``unnatural" complete graphs. For those cases, we can once again recall Cayley's formula: $\tau(K_n)=n^{n-2}$. Since the diameter of a cycle (resp. path) increases with the number of vertices, the condition $k < \text{diam } C_n = \lfloor n/2\rfloor$ is implicitly a lower bound on the number of vertices. Explicitly,
\begin{align*}
    \cali{C}_k&=\{{C_n}^k\ :\ n\ge 2k+1 \}\\  \cali{P}_k&=\{{P_n}^k\ :\ n\ge k+2\}.
\end{align*}
As discussed in the introduction, the initial values of a C-finite sequence determine the numerator of its rational generating function. The minimum values of $n$ above will mark the beginning of that family's sequence. Since paths and cycles are ubiquitous in graph theory, they are the focus of this paper. We provide a quick application below.


 Let $G$ be a connected graph. Suppose that we wish to estimate the number of spanning trees of its $k$-th power graph, $\tau(G^k)$. One way to estimate $\tau(G^k)$ is to fix a vertex $v\in V(G)$, let vertex sets $H_1,\dots,H_N$ be the components of $G-v$ and $m_i(v):= \max_{u\in H_i}|u,v|$. Then, deduce the bound:
$$\tau(G^k) \geq \prod_{i=1}^N \tau({P_{m_i(v)}}^k)$$
for all $v$. In particular,
$$\tau(G^k) \geq \max_{v\in V(G)} \prod_{i=1}^N \tau({P_{m_i(v)}}^k).$$


\section{Counting Total Number of Leaves} \label{sec:CountingLeaves}
Let $G$ be a graph and $T$ be a spanning tree of $G$. A \textit{leaf} of a tree $T$ is a vertex with degree exactly $1$ in $T$. We write $\cali{L}(T)$ for the set of leaves of $T$ and $\cali{T}(G)$ for the set of (labeled) spanning trees of $G$. When $G$ is clear from context, we simply write $\cali{T}$.

The next proposition plays a key role in our implementation for computing a parameter for graph families that we call the B-Z constant. The idea is to count the leaves (across all spanning trees) in a graph by counting how many times each vertex appears as a leaf. For each vertex in $G$, remove the vertex from $G$ and count the spanning trees in the resulting graph. (This method appears in \cite{aldous}.) 

\begin{prop}
    If $G$ is a labeled connected simple graph and for $v\in V(G)$, then
    $$\sum_{T\in \cali{T}} |\cali{L}(T)| = \sum_{v\in V(G)} \deg_G(v)\cdot |\cali{T}_v|$$
    where $\cali{T}:=\cali{T}(G)$ and $\cali{T}_v:=\cali{T}(G-v)$. 
\end{prop}
\begin{proof}\belowdisplayskip=-12pt 
     Fix $v\in V(G)$. Write $E_v:=\{e\in E(G)\ :\ v\in e\}$. There is a bijection between $E_v\times \cali{T}_v$ and the spanning trees of $T$ which contain $v$ as a leaf. Hence, we use indicators to obtain
     \begin{align*}
         \sum_{v\in V(G)} \deg_G(v)\cdot |\cali{T}_v| &=  \sum_{v\in V(G)} |\{T\in \cali{T}\ : v\in\cali{L}(T)|  \\
         &= \sum_{v\in V(G)}\sum_{T\in \cali{T}} \textbf{1}_{v\in \cali{L}(T)} = \sum_{T\in \cali{T}} |\cali{L}(T)|.
     \end{align*}\qedhere \end{proof}
We say $G$ is \textit{vertex-transitive} if for any $u,v\in V(G)$ there is an automorphism $\varphi$ of $G$ such that $\varphi(u)=v$ and $\varphi(v)=u$. 
\begin{cor}\label[cor]{cor:vtx-trans}
    If $G$ is vertex-transitive, then
    $$\sum_{T\in \cali{T}} |\cali{L}(T)| = n\cdot \deg_G(v)\cdot |\cali{T}_v|$$
    for any $v\in V(G)$.
\end{cor}

\section{The B-Z Constant: Average Number of Leaves per Vertex}\label{sec:bz-constant}

Next, we introduce a parameter for graph families whose member graphs are indexed by number of vertices. For such a graph family, our parameter represents the number of times (on average) that a vertex appears as a leaf in a spanning tree, averaged across all spanning trees, (asymptotically). 

For a graph family $\cali{G}$ indexed by number of vertices, where $G_n$ represents the graph with $n$ vertices in the family, we call the following the \textit{B-Z constant for} $\cali{G}$:
$$BZ(\cali{G}):=\lim_{n\to\infty}\frac{\sum_{T\in \cali{T}(G_n)} |\cali{L}(T)|}{n\cdot\tau(G_n)}$$
whenever the limit exists. This is the limit of the average number of leaves in a random spanning tree (and normalized by number of vertices) in the discussed family. For the infinite family of complete graphs $K_n$, the limit is famously equal to $1/e$ (see below).
In general, we have the bound $\frac{\sum |\cali{L}(T)|}{n\tau(G)} \le 1$ for any $G$. The rest of the section provides a few examples of B-Z constants and our approach for computing them.

The B-Z constants for both the family of cycles and the family of paths is equal to $0$, since every tree (a path) has a constant number of leaves (two). A star graph is a connected graph where all edges share the same vertex. When $\cali{G}$ is the family of star graphs, its B-Z constant is $1$ because the unique tree of a star graph (itself) has $n-1$ leaves.

As described in \Cref{sec:CountingLeaves}, it is possible to compute $\sum |\cali{L}(T)|$ by applying matrix-tree theorem. In the case of vertex-transitive graphs, computation of $\sum |\cali{L}(T)|$ (and hence the B-Z constant) is made faster thanks to \Cref{cor:vtx-trans}. We used this fact in the cases of: powers of cycles ($\cali{P}_k$) and torus graphs.

Incidentally, it's easy to compute the B-Z constant for complete graphs by using \Cref{cor:vtx-trans} in conjunction with Cayley's formula:

\begin{prop}
    The B-Z constant for the graph family $\{K_n\ :\ n\ge3\}$ is
    $\frac{1}{e}$.
\end{prop}
\begin{proof}
    Recall Cayley's formula, which states $\tau(K_n) = n^{n-2}$ for $n\ge 2$. With $\cali{T}:= \cali{T}(K_n)$, \Cref{cor:vtx-trans} tells us that $\sum |\cali{L}(T)| = n\cdot (n-1)^{n-2}$. Hence,
    $$\frac{\sum|\cali{L}(T)|}{n\cdot\tau(K_n)} = \left(\frac{n-1}{n}\right)^{n-2} = \left(1-\frac{1}{n}\right)^{n-2}$$
    which approaches $e^{-1}$ as $n\to\infty$.
\end{proof}

From a theorem in a joint paper by Hladký, Nachmias, and Tran \cite[Thm. 1.5]{leafdensity}, a stronger result follows. The B-Z constant for a family of graphs, where each member has minimum degree $\delta n$ (for a fixed $\delta > 0$), is $e^{-1}$. Their paper uses probabilistic methods, as well as graphons in the setting of dense graphs, to obtain this result. Each family $\cali{F}$ of graphs we considered in this paper is much sparser: there is a constant $c$ for which all $G\in \cali{F}$ satisfy $\Delta(G)\le c$; in other words, the maximum degree sequence $\{\Delta(G_n)\ : G_n\in \cali{F}\}$ is bounded by a constant. In the case of $\cali{G}_k$ and $\cali{P}_k$, a (tight) bound is $c= 2k$; whereas $c=4$ is tight for families of torus graphs and grid graphs.

Earlier, we noted that star graphs have B-Z constant equal to $1$. The next proposition observes that: for any rational number (in the interval $[0,1]$), we can find a graph family with B-Z constant equal to it. Such a family is constructed below by slightly modifying the family of star graphs.

\begin{prop}\label[prop]{prop:star-graph-bz-const}
    For any positive rational number $\frac{p}q < 1$, there is an indexed graph family $\cali{G}$ with $BZ(\cali{G})=p/q$.
\end{prop}
\begin{proof}
    First, we introduce a graph operation called subdivision. Let $e=uv$ be an edge of $G$, with $u,v\in V(G)$. Let $w\not\in V(G)$. The result of subdividing $e$ in $G$ is the graph $G'$ where $V(G') := V(G) \cup \{w\}$ and $E(G') := (E(G)\setminus \{e\})\cup \{uw,wv\}$.

    To construct each member $G_n\in\cali{G}$ (assuming $n\ge 1$): begin with a star graph $S_{pn}$ with $pn$ leaves; then, subdivide $(q-p)n$ edges in $S_{pn}$ and call this $G_n$. Note that $|V(G_n)| = qn+1$ and that $G_n$ has $pn$ leaves. Furthermore, $G_n$ is a tree so it has exactly one spanning tree. It follows that the B-Z constant for $\cali{G}$ is $\lim_{n\to\infty} \frac{pn}{qn+1} = \frac{p}q$.\end{proof}

Unlike these modified star graph families, the families we consider in upcoming sections have B-Z constants that are algebraic, but not rational.

\section{Powers of a Cycle: Experimental Results} \label{sec:cycles}


  In upcoming sections, we determine the B-Z constants of several graph families. It's possible, as we do in later sections, to compute the B-Z constant of a graph family from its asymptotic behavior in the total number of leaves and number of spanning trees (which can be obtained from their respective generating functions). We used standard residue calculations, implemented in procedure \T{BZc} from our Maple package (\Cref{maple-package}).


\subsection{Generating Functions for the Number of Spanning Trees}\label{cycle-tree-gf}


  In upcoming sections, we determine the B-Z constants of several graph families. It's possible, as we do in later sections, to compute the B-Z constant of a graph family from its asymptotic behavior in the total number of leaves and number of spanning trees (which can be obtained from their respective generating functions). We used standard residue calculations, implemented in procedure \T{BZc} from our Maple package (\Cref{maple-package}).


\subsection{Generating functions for the Total Number of Leaves} \label{cycle-leaf-gf}
\begin{thm}
    The generating function for the total number of leaves (across all spanning trees of a member) in $\cali{C}_2$ is
    $$\frac{-8(10t^7 - 67t^6 + 109t^5 + 99t^4 - 282t^3 - 30t^2 + 145t - 40)}{(t + 1)^2(t^2 - 3t + 1)^3}.$$
\end{thm}

\begin{thm}
The generating function for the total number of leaves (across all spanning trees of a member) in $\cali{C}_3$ is
    $$\frac{A_3}{(t - 1)^3(t^4 + 3t^3 + 6t^2 + 3t + 1)^3(t^4 - 4t^3 - t^2 - 4t + 1)^3}$$
where 
\begin{align*}
    A_3&:=-8820t^{26} + 51390t^{25} + 61812t^{24} + 2088t^{23} - 2539950t^{22} - 2981160t^{21} + 2492784t^{20} \\ 
    &+ 45845688t^{19}+ 83018808t^{18}+ 107694630t^{17} - 44892840t^{16} - 166389300t^{15} - 333210654t^{14}\\ & - 121438506t^{13} +42702660t^{12}+ 312824052t^{11}+ 213402930t^{10} + 100784592t^9 - 77616756t^8\\ & - 90041700t^7 - 62209728t^6- 13836186t^5+ 276924t^4 + 2761596t^3+ 501534t^2 + 32592t - 54432.
\end{align*}
\end{thm}
As the reader might have noticed, the numerators for these generating functions become too cumbersome to write (and more quickly than in the previous section). Henceforth, we omit the numerator and refer the reader to this website \cite{spt} for detailed results when $4\leq r\leq 5$.

\begin{thm}
The generating function for the total number of leaves (across all spanning trees of a member) in $\cali{C}_4$ is
    $$\frac{A_4}{B_4}$$
where $A_4$ is a polynomial of degree $80$ and 
\begin{align*}
    B_4&=\\ &(t + 1)^3(t^6 - 3t^5 + 6t^4 - 10t^3 + 6t^2 - 3t + 1)^3(t^8 - 4t^7 - 17t^6 + 8t^5 + 49t^4 + 8t^3 - 17t^2 - 4t + 1)^3\\
    &(t^{12} + 3t^{11} + 12t^{10} + 28t^9 - 27t^8 + 36t^7 - 81t^6 + 36t^5 - 27t^4 + 28t^3 + 12t^2 + 3t + 1)^3
\end{align*}

\end{thm}

\subsection{B-Z Constants} \label{BZ-cycle}
See \Cref{sec:bz-constant} to read about our method for computing these constants. Recall that $BZ(\cali{G})$ denotes the B-Z constant of a graph family $\cali{G}$.

\begin{thm}
    The B-Z constant for $\cali{C}_2$ is 
    $$-6+\frac{14\sqrt{5}}{5}.$$
\end{thm}
\begin{thm}
    The B-Z constant for $\cali{C}_3$ is 
    $$\frac37 \left(\frac{45}2 + 9 \sqrt{7} - \frac12 \sqrt{3857 + 1684 \sqrt{7}}\right).$$
\end{thm}

\begin{thm}
    Let $\alpha$ be the smallest real root of $z^8-4z^7-17z^6+8z^5+49z^4+8z^3-17z^2-4z+1$. Then, the B-Z constant for $\cali{C}_4$ is 
   $$\frac{1}{2025}\left(-216\alpha^7-2144\alpha^6+16344\alpha^5+41056\alpha^4-17064\alpha^3-87936\alpha^2-25304\alpha+6944\right).$$
\end{thm}
\begin{cor}
    $$BZ(\cali{C}_3) < \frac{117451}{355635} < BZ(\cali{C}_4) < \frac{117452}{355635}$$
\end{cor}
In the previous theorem, $\alpha \approx 0.158778$.

\section{Powers of a Path: Experimental Results} \label{sec:paths}
\subsection{Generating functions for the Number of Spanning Trees} \label{path-tree-gf}

In this section, we present our results for the Number of Spanning Trees in the class $\cali{P}_k$, with $2\le k \le 6$.

\begin{thm}
    The generating function $f(t)$ for the number of spanning trees in $\cali{P}_2$ is
    $$\frac{-3t+8}{t^{2} - 3t + 1}.$$
\end{thm}
   
\begin{thm}
    The generating function $f(t)$ for the number of spanning trees in $\cali{P}_3$ is
   $$\frac{-16t^4 + 77t^3 - 33t^2 + 39t - 75}{(t - 1)(t^{4} - 4t^{3} - t^{2} - 4t + 1)}.$$
\end{thm}

\begin{thm}
    The generating function $f(t)$ for the number of spanning trees in $\cali{P}_4$ is
   $$\frac{M_4}{(t^{6} - 3t^{5} + 6t^{4} - 10t^{3} + 6t^{2} - 3t + 1)(t^{8} - 4t^{7} - 17t^{6} + 8t^{5} + 49t^{4} + 8t^{3} - 17t^{2} - 4t + 1)}.$$
   where 
   \begin{align*}
       M_4&=-125t^{13} + 859t^{12} - 13t^{11} - 3141t^{10} + 3475t^9 - 5968t^8 - 11312t^7 \\
       &+ 36080t^6 - 5597t^5 - 7893t^4 + 2435t^3 - 2741t^2 - 413t + 864\\
   \end{align*}
\end{thm}

We know the generating functions for larger $r$ values, but there are many terms with longer coefficients. So, we do not state them in the paper. They are stated in this website \cite{spt}.

\begin{thm}
    The generating function $f(t)$ for the number of spanning trees in $\cali{P}_5$ is
   $$\frac{M_5}{E_5}.$$
   where 
   $M_5$ is a degree $40$ polynomial in $t$ and
   \begin{align*}
       E_5&=(t - 1)(t^{8} + 3t^{7} + 6t^{6} - t^{5} + 15t^{4} - t^{3} + 6t^{2} + 3t + 1)(t^{16} - 5t^{15} + 10t^{14} - 10t^{13} - 28t^{12}\\
       & + 10t^{11}+ 110t^{10} + 110t^{9} + 88t^{8} + 110t^{7} + 110t^{6} + 10t^{5} - 28t^{4} - 10t^{3} + 10t^{2} - 5t + 1)(t^{16}\\
       &- 5t^{15} - 23t^{14}- 10t^{13} - 94t^{12} - 485t^{11} + 242t^{10} + 110t^{9} + 649t^{8} + 110t^{7} + 242t^{6} - 485t^{5}\\
       &- 94t^{4} - 10t^{3} - 23t^{2}- 5t + 1)    
   \end{align*}
\end{thm}


\subsection{Generating functions for the Total Number of Leaves} \label{path-leaf-gf}

\begin{thm}
The generating function for the total number of leaves (across all spanning trees of a member) in $\cali{P}_2$ is
    $$\frac{-2(2t^3 - 15t^2 + 27t - 9)}{(t^2 - 3t + 1)^2}$$
\end{thm}

\begin{thm}
The generating function for the total number of leaves in $\cali{P}_3$ is
    $$\frac{2(16t^{10}-154t^9+403t^8-340t^7+963t^6-768t^5+1109t^4-788t^3+509t^2-470
t+96)}{(t-1)^2(t^4-4t^3-t^2-4t+1)^2}$$
\end{thm}

\begin{thm}
The generating function for the total number of leaves in $\cali{P}_4$ is
    $$\frac{M_4}{(t^6 - 3t^5 + 6t^4 - 10t^3 + 6t^2 - 3t + 1)^2(t^8 - 4t^7 - 17t^6 + 8t^5 + 49t^4 + 8t^3 - 17t^2 - 4t + 1)^2}$$
    where $M_4$ is a degree $29$ polynomial in $t$. Note that the degree of the denominator is $28.$ 
\end{thm}

\subsection{B-Z Constants} \label{BZ-path}
See \Cref{sec:bz-constant} to read about our method for computing these constants. Compare the results in this section to \Cref{BZ-cycle}.

\begin{thm}
    The B-Z constant for $\cali{P}_2$ is 
    $$-6+\frac{14\sqrt{5}}{5}.$$
\end{thm}
\begin{thm}
    The B-Z constant for $\cali{P}_3$ is 
    $$\frac37 \left(\frac{45}2 + 9 \sqrt{7} - \frac12 \sqrt{3857 + 1684 \sqrt{7}}\right).$$
\end{thm}

\begin{thm}
    $BZ(\cali{P}_4)=BZ(\cali{C}_4)$.
\end{thm}

We also verified that $BZ(\cali{P}_5)=BZ(\cali{C}_5)$. See \cite{spt}.

\section{B-Z Constants for Grid and Torus Graphs}\label{sec:gridtorus}
In previous sections, we discussed powers of cycles and powers of paths. The members of these two graph families were related by a subgraph relation and we found their B-Z constants to be the same up to $r=5$. Since a power of a cycle is vertex-transitive, computing the total number of leaves was faster (\Cref{cor:vtx-trans}). In this section, we discuss another pair of graph families -- grid graphs and torus graphs -- which are also related by a subgraph relation. Similarly to the previous pair of families, torus graphs are vertex-transitive.  However, the B-Z constants for grid graphs and torus graphs differ.

\textbf{Definition.} The $a\times b$ grid graph $P_a\times P_b$ has vertex set $V(P_a\times P_b):=\{(i,j)\ :\ i\in[a], j\in[b]\}$ and two (distinct) vertices $(u_1,u_2)$ and $(v_1,v_2)$ are adjacent whenever $|u_1-v_1|=1$ or $|u_2-v_2|=1$, but not both.

\textbf{Definition.} The $a\times b$ torus graph $C_a\times C_b$ has vertex set $V(C_a\times C_b):=\{(i,j)\ :\ i\in[a], j\in[b]\}$ and two (distinct) vertices $(u_1,u_2)$ and $(v_1,v_2)$ are adjacent whenever $|u_1-v_1|\in\{1,a-1\}$ or $|u_2-v_2|\in\{1,b-1\}$, but not both.

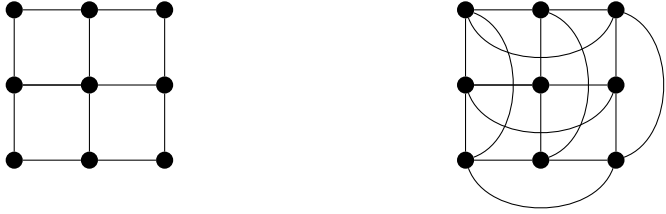
\begin{figure}[H] 
   \centering
    \begin{tikzpicture}

\tikzstyle{knode}=[circle,fill,draw=black,thick,inner sep=2pt] 

\pgfmathsetmacro{\n}{7} 
\pgfmathsetmacro{\R}{1cm} 
\pgfmathsetmacro{\A}{90} 

\node (g1) at (-1,-1) [knode] {};
\node (g2) at (0,-1) [knode] {};
\node (g3) at (1,-1) [knode] {};
\node (g4) at (-1,0) [knode] {};
\node (g5) at (0,0) [knode] {};
\node (g6) at (1,0) [knode] {};
\node (g7) at (-1,1) [knode] {};
\node (g8) at (0,1) [knode] {};
\node (g9) at (1,1) [knode] {};

\draw (g1)--(g2);
\draw (g3)--(g2);
\draw (g1)--(g4);
\draw (g7)--(g4);
\draw (g5)--(g4);
\draw (g5)--(g4);
\draw (g7)--(g8);
\draw (g9)--(g8);
\draw (g5)--(g6);
\draw (g5)--(g8);
\draw (g9)--(g6);
\draw (g2)--(g5);
\draw (g3)--(g6);

\node (t1) at (-1+6,-1) [knode] {};
\node (t2) at (0+6,-1) [knode] {};
\node (t3) at (1+6,-1) [knode] {};
\node (t4) at (-1+6,0) [knode] {};
\node (t5) at (0+6,0) [knode] {};
\node (t6) at (1+6,0) [knode] {};
\node (t7) at (-1+6,1) [knode] {};
\node (t8) at (0+6,1) [knode] {};
\node (t9) at (1+6,1) [knode] {};

\draw (t1)--(t2);
\draw (t3)--(t2);
\draw (t1)--(t4);
\draw (t7)--(t4);
\draw (t5)--(t4);
\draw (t5)--(t4);
\draw (t7)--(t8);
\draw (t9)--(t8);
\draw (t5)--(t6);
\draw (t5)--(t8);
\draw (t9)--(t6);
\draw (t2)--(t5);
\draw (t3)--(t6);
\draw (t1) to[out=-70,in=-110] (t3);
\draw (t4) to[in=-110,out=-70] (t6);
\draw (t7) to[in=-110,out=-70] (t9);

\draw (t9) to[out=-20,in=20] (t3);
\draw (t8) to[out=-20,in=20] (t2);
\draw (t7) to[out=-20,in=20] (t1);

\end{tikzpicture}
    \caption{On the left, the $3\times 3$ grid graph. On the right, the $3\times 3$ torus graph.}
    \label{fig:fig2}
\end{figure}

For the family of torus graphs, we assume $\min(a,b) \ge 3$ to avoid considering multigraphs. A torus graph (resp. grid graph) with dimensions $a\times b$ is isomorphic to another torus graph (resp. grid graph) with dimensions $b\times a$; hence, we may treat the two interchangeably. 

\begin{thm}
    The B-Z constant for the family of $2\times n$ grid graphs is
    $$\sqrt{3}-\frac{3}{2}.$$
\end{thm}

In the $3\times n$ case, the B-Z constants of the two corresponding families differ.   

\begin{thm} \label{torus3}
    The B-Z constant for the family of $3\times n$ torus graphs is
    $$\frac{24\sqrt{21}-104}{21}.$$
\end{thm}

\begin{thm} \label{grid3}
    The B-Z constant for the family of $3\times n$ grid graphs is
    $$\frac{-10500 + 2200 \sqrt{21} + 1197 \sqrt{230 - 50 \sqrt{21}} + 149 \sqrt{210 (23 - 5 \sqrt{21})}}{5040}.$$
\end{thm}

Similarly, in the case of $4\times n$ grid and torus graph, their B-Z constants differ again \cite{spt}. We conjecture that they differ for all $r\times n$ when $r\ge 3$.  

The algebraic numbers of the families above increase in complexity as the fixed dimension of the corresponding family increases. In this website \cite{spt}, we recorded the outputs of these results -- symbolically and numerically. We summarize some of those results here, with decimal approximations (and at least six digits):

The B-Z constant for: $4\times n$ torus graphs is $0.2917148\dots$; $5\times n$ torus graphs is $0.29342497\dots$; and $6\times n$ torus graphs is $0.294011\dots$. Notably, these B-Z constants appear to approach the B-Z constant for $n\times n$ torus graphs -- which is remarked to be $(1-\frac2\pi)\cdot\frac{8}{\pi^2}$ by \cite[Remark 5.3]{leafdensity}, who furthermore cite \cite[p. 112]{lyonsperes}. The B-Z constant for $4\times n$ grid graphs is $.2746417\dots$.

We computed B-Z constants for fewer grid families (than torus families) because grid graphs are not vertex-transitive; hence, we were unable to apply \Cref{cor:vtx-trans} to speed up computation of the ``total number of leaves".

\section{Estimating B-Z Constants with Random Sampling}\label{sec:verification}

As we have shown in previous sections, with enough computational resources, it is possible to find the exact B-Z constant for several graph families. Even though the exact B-Z constant may be given as a root of a polynomial, one can also obtain precise numerical bounds for it. To perform these exact computations, we wrote a procedure to compute the total number of leaves of a graph (\texttt{NumLeaves} or \texttt{VtxTransNumLeaves}) and another to compute the number of spanning trees (\texttt{NumSpanTree}). The ratio between the outputs of these procedures gave us the average number of leaves of a graph. 

However, if one only cares about finding a numerical estimate of the B-Z constant for a graph family -- after confirming it exists -- there is another approach. Fix a member of the graph family with a large number $n > N$ of vertices. Then, sample a large number $K$ of (uniformly) random spanning trees from the graph $G_n$, count the number of leaves in each sampled tree, and use the resulting data to estimate the average number of leaves for the large member $G_n$ in the family. The B-Z constant is an asymptotic value; so, for large enough $N$ and $K$, a numerical sample will suffice. 

In our Maple package \cite{spt}, we implemented a procedure (\T{EstAvgNumLeaves}) that uses parameters $N,K$, as described above, to estimate the average number of leaves in such a graph $G_n$. Dividing the output from \T{EstAvgNumLeaves} by the number of vertices ($n$) will result in a numerical estimate for the B-Z constant. To sample uniformly random spanning trees, we implemented Wilson's algorithm, described in \cite{wilsonalg}.

Using the procedure, we compared numerical estimates to the exact B-Z constants for the graph families in this paper. We used $N\ge 100$ and $K\ge 100$ (and usually $K\ge 200$). As expected, we found that the estimated and the true B-Z constants were similar. The data from these numerical estimates can be found in \cite{spt}. 

\section{Accompanying Maple package} \label{maple-package}
The \T{Help()} procedure in the Maple package, available at \cite{spt}, will provide even more details, as well as examples for how to use these procedures.

In the following procedures, graphs are represented as an \T{exprseq} type in Maple. Particularly, a graph is \T{n,E}, where \T{n} is a positive integer and \T{E} is a set of edges on $[n]$. We assume all graphs are connected.

Our Maple package broadly depends on the \texttt{LinearAlgebra} library. The \texttt{RandomTree} procedure also makes use of the \texttt{GraphTheory} library. Below we list some key procedures along with their descriptions:

\texttt{NumSpanTree(n,E)} given a positive integer \texttt{n} and a set of edges \T{E} on the set $\{1,\dots,n\}$, the procedure returns the number of spanning trees of the corresponding graph $n,E$.

\texttt{NumSpanTreeSeq(F, ArgList, a,b)} given the name of a graph-generating procedure \T{F} and a list of arguments \T{ArgList}, outputs a list of the number of spanning trees for the graphs \T{F(a, op(ArgList))},$\dots$,\T{F(b, op(ArgList))}. 


\T{NumLeaves(n,E)} given a graph \T{n,E}, returns the total number of leaves in such a graph across all its spanning trees.

\T{NumLeavesSeq(F, ArgList, a,b)} given the name of a graph-generating procedure \T{F} and a list of arguments \T{ArgList}, outputs a list of the total number of leaves for the graphs \T{F(a, op(ArgList))},$\dots$,\T{F(b, op(ArgList))}. 

\T{VtxTransNumLeavesSeq(F, ArgList, a,b)} similar to \T{NumLeavesSeq}. However, assumes that the graphs generated by \T{F} are vertex-transitive. Uses \Cref{cor:vtx-trans} to compute the output faster.


\T{Hnr(n,r)} constructs the $r$-th power of a path on $n$ vertices. Returns $n,E$.

\T{Gnr(n,r)} constructs the $r$-th power of a cycle on $n$ vertices. Returns $n,E$.


\T{BZc(F,Arglist,a,b,k)} Using a graph-generating procedure \T{F} with fixed arguments \T{ArgList} and varying the first argument from \T{a} to \T{b}, compute the B-Z constant. \T{k} adjusts the B-Z computation so that if the index of a graph (in the family) is $n$, then the graph has $k\cdot n$ vertices.

The following procedures generate the outputs for the (main) theorems in this paper, up to \T{R} $\ge 2$, with \T{K} terms in the sequence, and in variable \T{x}: 

\T{HnrS(R,x,K)}

\T{GnrS(R,x,K)}

From Doron Zeilberger's \T{Cfinite.txt} package: 

\T{RGF(S,t)} given a C-finite sequence \T{S}, outputs its rational generating function in the variable \T{t}.

\section{Conjectures and Further Discussion} \label{sec:conj}
Some observations from our results may indicate deeper structural insights:
\begin{enumerate}
    \item For each $r$ we computed, the denominator of a rational function from \Cref{path-tree-gf} divides the (corresponding) denominator of a rational function from \Cref{cycle-tree-gf}. Similarly for \Cref{path-leaf-gf} and \Cref{cycle-leaf-gf}.  
        \item In the denominators for the generating functions given in \Cref{cycle-tree-gf}, \Cref{cycle-leaf-gf}, \Cref{path-tree-gf}, \Cref{path-leaf-gf}, each polynomial factor has coefficients which are symmetric. Is this also true for all other recursive graph families?
\end{enumerate}

While in this paper we mainly provided results for powers of cycles and powers of paths, we also computed the corresponding rational generating functions for grid graphs and torus (grid) graphs (see \cite{spt}). Those two classes of families, when compared to each other, did not satisfy property 1.; however, they did satisfy property 2 (as listed above). Since a grid graph is a subgraph of a torus graph (of the same dimensions), we know that property 1. is not satisfied from just a subgraph relation. We conclude with a few conjectures. 

\begin{conj}
    The B-Z constants of $\cali{C}_k$ and $\cali{P}_k$ are equal for all $k\ge 1$.
\end{conj}

Assuming that the previous conjecture holds, a natural question arises: given a graph family, is there a threshold for how many vertices you can "modify" (by edge deletion) while keeping the same B-Z constant? It turns out that grid graphs and torus graphs have different B-Z constants, which may suggest that we can only modify a finite number of vertices. Originally, we considered the following conjecture. Note that the conditions in the following conjecture are satisfied when $\cali{G}=\cali{C}_k$ and $\cali{H}=\cali{P}_k$ (as defined in \Cref{graph-fams}).

\begin{conj}
    \textbf{\emph{(False)}} Let $\cali{G}:=\{G_n\}_n$ and $\cali{H}:=\{H_n\}_n$ be (indexed) families of connected graphs, such that $BZ(\cali{G})$ and $BZ(\cali{H})$ both exist. If the following hold:
    \begin{itemize}
        \item $H_n\subseteq G_n$ for all $n$ and
        \item there is some $c_1$ for which  $|E(G_n\setminus H_n)|\leq c_1$ for all $n$,
    \end{itemize}
   then $\cali{G}$ and $\cali{H}$ have the same B-Z constant. 
\end{conj}

For two graphs $G,H$ such that $V(G)=V(H)$, we define $G\setminus H$ as the graph $V(G\setminus H):=V(G)$ and $E(G\setminus H):=E(G)\setminus E(H)$. Unfortunately, the conjecture above fails for families of sparse graphs, as proven by considering the following counterexample.

\textbf{Counterexample.} Let $\cali{H}$ be the family of graphs constructed in the proof of \Cref{prop:star-graph-bz-const} with B-Z constant equal to $\frac12$. In other words, $\cali{H}$ is the family of graphs resulting from subdividing each edge in a star graph once. From each $H_n\in \cali{H}$, construct $G_n$ by adding an edge between two non-leaf vertices. We see that $\tau(G_n)=3\cdot\tau(H_n)$, but $\cali{L}(G_n)=\cali{L}(H_n)$. Hence, $BZ(\cali{G})=\frac{1}{3}\cdot BZ(\cali{H}) = \frac16 \neq \frac12 = BZ(\cali{H})$.

\begin{obs}
    Let $n\ge 2r+1$. If $G_n\in \cali{C}_k$ and $H_n\in \cali{P}_k$, then $G_n\setminus H_n$ contains a (copy of the) half graph on $2r$ vertices and $n-2r$ isolated vertices.
\end{obs}
The \textit{half graph} on $2n$ vertices is a bipartite graph with vertex set $\{u_1,\dots,u_n,v_1,\dots,v_n\}$ and edges $\{u_iv_j\ :\ i\leq j\}$. Since $r$ is fixed, the number of differing edges between a power of a cycle and a power of a path is finite.

\typeout{get arXiv to do 4 passes: Label(s) may have changed. Rerun}

\begin{center}{\textbf{Acknowledgments}}
\end{center}
The authors wish to thank an anonymous reviewer for their detailed and insightful comments that improved the readability of this paper, as well as for providing additional references to the paper (which relate it to pre-existing results and literature).

\bibliographystyle{abbrv}

\end{document}